\documentclass[reqno]{amsart}

\usepackage{graphicx}
\usepackage{amscd}
\usepackage{amsmath, amssymb}
\usepackage{graphics}
\usepackage{graphicx}
\usepackage{epsfig}
\usepackage{xcolor}
\usepackage{soul}
\usepackage{tikz}

\newtheorem{theorem}{Theorem}

\newtheorem{definition}[theorem]{Definition}

\newtheorem{proposition}[theorem]{Proposition}

\begin{document}

\title[High-dimensional representations of $T_n$ via several techniques]{Constructing high-dimensional representations of twin groups via several techniques}

\author{Mohamad N. Nasser}         
\address{Mohamad N. Nasser\\
         Department of Mathematics and Computer Science\\
         Beirut Arab University\\
         P.O. Box 11-5020, Beirut, Lebanon}
         \email{m.nasser@bau.edu.lb}

\maketitle
 
\begin{abstract}
In this article, we consider the twin group $T_n$ for a positive integer $n$ and present two techniques for constructing high-dimensional and complicated representations of $T_n$. The first technique consists of composing a known representation of $T_n$ with the well-established Cohen map, while the second technique relies on the tensor product of two known representations of $T_n$. To implement these techniques, we consider the two representations $\eta_1$ and $\eta_2$ of $T_n$ introduced recently by Nasser. The first technique is applied to the representation $\eta_1: T_n \longrightarrow \mathrm{Aut}(\mathbb{F}_n)$, which maps $T_n$ into the automorphism group of the free group $\mathbb{F}_n$ with $n$ generators. The second technique is applied to the representation $\eta_2: T_n \longrightarrow \mathrm{GL}_n(\mathbb{Z}[t^{\pm 1}])$, which maps $T_n$ into the general linear group over the ring of Laurent polynomials in an indeterminate $t$. Using these techniques, we construct new high-degree representations of $T_n$ and examine, for some of them, structural properties such as faithfulness and reducibility. This work provides new avenues for developing representations of $T_n$ and related groups using these techniques.
\end{abstract}

\renewcommand{\thefootnote}{}

\vspace{0.5cm}

\section{Introduction}

\begin{figure}[h!] 
\centering
\begin{minipage}{0.45\textwidth}
\centering
\begin{tikzpicture}[scale=0.9]
	\draw[thick] (-1.5,0)--(-1.5,2); 
    \fill (-1,1) circle(1pt) (-1.2,1) circle(1pt)(-0.8,1)circle(1pt);       
    \draw[thick] (-.5,0)--(-.5,2);
    \draw[thick] (2.5,0)--(2.5,2); 
    \fill (2,1) circle(1pt) (2.2,1) circle(1pt) (1.8,1)circle(1pt);       
    \draw[thick] (1.5,0)--(1.5,2);
    \draw[thick] (0,0) to[out=90, in=-90] (1,2);
    \draw[thick, white, line width=4pt] (1,0) to[out=90, in=-90] (0,2);
    \draw[thick] (1,0) to[out=90, in=-90] (0,2);
    \node[above] at(0,2){$i$};
    \node[above] at(1,2){$i+1$};
\end{tikzpicture}

\vspace{0.1cm}
{\small The generator $\sigma_i$}
\end{minipage}
\hfill
\begin{minipage}{0.45\textwidth}
\centering
\begin{tikzpicture}[scale=0.9]
\draw[thick] (5.5,0)--(5.5,2); 
\fill (5,1) circle(1pt) (4.8,1) circle(1pt)(5.2,1)circle(1pt);       
\draw[thick] (4.5,0)--(4.5,2);
\draw[thick] (8.5,0)--(8.5,2); 
\fill (8,1) circle(1pt) (8.2,1) circle(1pt)(7.8,1)circle(1pt);       
\draw[thick] (7.5,0)--(7.5,2);    
\draw[thick] (7,0) to[out=90, in=-90] (6,2);
\draw[thick] (6,0) to[out=90, in=-90] (7,2);
\node[above] at(6,2){$i$};
\node[above] at (7,2) {$i+1$};
\end{tikzpicture}

\vspace{0.1cm}
{\small The generator $s_i$}
\end{minipage}
\caption{The braid generator $\sigma_i$ and the twin generator $s_i$.}
\label{Fig1}
\end{figure}
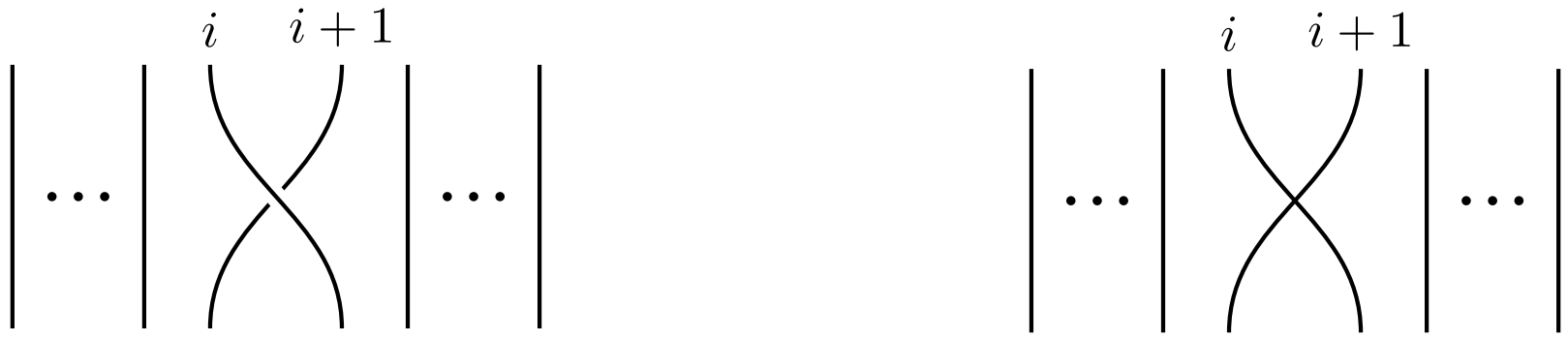

The twin group on $n$ strands, denoted by $T_n$, is a Coxeter group generated by the $n-1$ involutions $s_1,s_2,\ldots,s_{n-1}$ satisfying the defining relations
\begin{align*} 
s_i^2 & =1, \quad \ \qquad i=1,2,\ldots,n-1,\\ 
s_is_j& =s_js_i, \hspace*{0.02cm} \qquad |i-j|\ge 2. 
\end{align*}
This group was first introduced by G. Shabat and V. Voevodsky in \cite{Shabat1990} and has also been studied in the literature under several different names, such as the flat braid group and the planar braid group; see, for example, \cite{Khovanov1996, Khovanov1997, Merkov1999}. The group $T_n$ admits a geometric interpretation analogous to that of the classical braid group on $n$ strands, denoted by $B_n$. Despite this similarity, the two groups differ significantly in their structure. In $T_n$, the generators $s_i$ represent simple transposition-type interactions between strands $i$ and $i+1$ without distinguishing between over- and under-crossings. On the other hand, the generators $\sigma_i$ of $B_n$ correspond to genuine crossings in which strand $i$ passes over strand $i+1$, and they satisfy the braid relations $\sigma_i\sigma_{i+1}\sigma_i=\sigma_{i+1}\sigma_i\sigma_{i+1}, 1\le i\le n-2,$ which do not hold in $T_n$. Therefore, although both groups arise from configurations of $n$ strands, their algebraic structures and geometric interpretations are fundamentally different (see Figure \ref{Fig1}).

Representation theory plays a fundamental role in the study of algebraic structures by realizing their elements as automorphisms of free groups or as matrices acting on vector spaces. Such realizations make it possible to analyze the properties of these structures using tools from linear algebra and often reveal structural features that are not easily accessible through purely combinatorial approaches. Specifically, the construction of higher-dimensional or more complicated representations is important for gaining deeper insight into the algebraic structure. For instance, representations of greater complexity have a higher likelihood of being faithful and irreducible, thereby providing more precise information about the algebraic structure. \vspace*{0.15cm}

In this paper, we introduce two techniques to construct higher-dimensional and more complicated representations of the twin group $T_n$. The first technique involves composing an existing $n$-dimensional representation of $T_n$ with the well-known Cohen map described in \cite{Cohen1996}, resulting in a new representation of dimension $nk$, where $k$ is a positive integer. This technique has been applied on many representations of the braid group $B_n$ and its pure subgroup $P_n$ see, for instance, \cite{Abdulrahim2003, Abdulrahim2005,Nasser2020}. The second technique utilizes the tensor product of representations: given two matrix representations of $T_n$ with degrees $n_1$ and $n_2$, we form their tensor product by taking the Kronecker product of the matrices. This produces a new representation acting on the tensor product of the corresponding vector spaces, yielding a matrix representation of $T_n$ of degree $n_1 n_2$. \vspace*{0.15cm}

In our research, we consider two well-known representations of the twin group $T_n$ described explicitly in \cite{Nasser2025}, with the aim of constructing higher-dimensional and more complicated representations of $T_n$ from them using the above two techniques. The first representation, $\eta_1: T_n \longrightarrow \text{Aut}(\mathbb{F}_n)$, maps $T_n$ into the automorphism group of the free group $\mathbb{F}_n$ with $n$ generators. The second representation, $\eta_2: T_n \longrightarrow \text{GL}_n(\mathbb{Z}[t^{\pm 1}])$, maps $T_n$ into the general linear group over the ring of Laurent polynomials in an indeterminate $t$. We recall the definitions and main properties of these representations in Section 2. \vspace*{0.15cm}

Our main results are presented separately in Sections 3 and 4. First, in Section 3, we describe in detail the three steps required to implement the first technique, and we apply these steps to the representation $\eta_1$ to construct a new representation $\widetilde{\eta}_1: T_n \longrightarrow \mathrm{GL}_{2n}(\mathbb{Z}[t^{\pm1}])$, where $t$ is an indeterminate (Proposition \ref{propppp} and Theorem \ref{Th1st}). In Section 4, we apply the second technique to the representation $\eta_2$ of $T_n$ for the cases $n=2$ and $n=3$, obtaining two representations $\hat{\eta}_2$ and $\bar{\eta}_2$ of dimensions $4$ and $9$, respectively (Propositions \ref{firstprop} and \ref{secondprop}). Furthermore, we examine the faithfulness and reducibility of these two representations.

\section{Twin Groups Representations}

In this section, we recall the definition and main properties of two representations of the twin group $T_n$ constructed in \cite{Nasser2025}. We begin with the representation $\eta_1$, which is defined from $T_n$ to $\mathrm{Aut}(\mathbb{F}_n)$, where $\mathbb{F}_n$ denotes the free group on $n$ generators. We also present the matrix representation derived from it.

\begin{definition} \cite{Nasser2025} \label{defeta1}
The representation  $\eta_1: T_n \longrightarrow \mathrm{Aut}(\mathbb{F}_n)$, where $\mathbb{F}_n$ is a free group with $n$ generators $x_1,x_2,\ldots, x_n$, is the representation that acts on the generators $s_i, 1\leq i \leq n-1,$ as follows.
$$\eta_1(s_i):
\left\{\begin{array}{l}
x_i\rightarrow x_ix_{i+1}x_i^{-1},\\
x_{i+1}\rightarrow x_ix_{i+1}^{-1}x_ix_{i+1}x_i^{-1},\\
x_j\rightarrow x_j \ \text{ for } \ j \notin \{i,i+1\} . \\
\end{array}\right.$$
\end{definition}

\begin{theorem} \cite{Nasser2025} \label{defeta1matrix}
The matrix representation of \(T_n\) derived from \(\eta_1\), via Magnus representation and Fox derivatives, is the representation $\eta_1: T_n \longrightarrow \mathrm{GL}_n(\mathbb{Z}[t^{\pm 1}]),$ where \(t\) is an indeterminate. For simplicity, we still denote it by \(\eta_1\). The action of $\eta_1$ on the generators $s_i, 1 \le i \le n-1$, is given as follows.
$$\eta_1(s_i)= \left( \begin{array}{c|@{}c|c@{}}
   \begin{matrix}
     I_{i-1} 
   \end{matrix} 
      & 0 & 0 \\
      \hline
    0 &\hspace{0.2cm} \begin{matrix}
   	1-t & t\\
   	2-t & t-1\\
\end{matrix}  & 0  \\
\hline
0 & 0 & I_{n-i-1}
\end{array} \right) \text{ for } 1\leq i \leq n-1.$$ 
\end{theorem}

\begin{theorem}\cite{Nasser2025}
The representation $\eta_1: T_n \longrightarrow \mathrm{GL}_n(\mathbb{Z}[t^{\pm 1}])$ is faithful for $n=2,3$.
\end{theorem}
 
\begin{theorem} \cite{Nasser2025}
The representation $\eta_1: T_n \longrightarrow \mathrm{GL}_n(\mathbb{Z}[t^{\pm 1}])$ is reducible to the degree $n-1$ for all $n\geq 3$. Moreover, the complex specialization of its $(n-1)$-composition factor, namely $\eta_1': T_n \longrightarrow \mathrm{GL}_{n-1}(\mathbb{C})$, where $t\neq 0$, is irreducible if and only if $t\neq \dfrac{2n-2}{n-2}$ and $t\neq 2$.
\end{theorem}

We now recall the representation $\eta_2$ of $T_n$, which forms a family of representations of $T_n$ to $\mathrm{GL}_n(\mathbb{Z}[t^{\pm 1}])$, where $t$ is an indeterminate. This family depends on an invertible element $f(t) \in \mathbb{Z}[t^{\pm 1}]$.

\begin{definition}\cite{Nasser2025} \label{defeta2}
The representation $\eta_2: T_n \longrightarrow \mathrm{GL}_n(\mathbb{Z}[t^{\pm 1}])$, where $t$ is an indeterminate, is the representation that acts on the generators $s_i, 1\leq i \leq n-1,$ as follows.
$$\eta_2(s_i) = \left( \begin{array}{c|@{}c|c@{}}
   \begin{matrix}
     I_{i-1} 
   \end{matrix} 
      & 0 & 0 \\
      \hline
    0 &\hspace{0.2cm} \begin{matrix}
   	0 & f(t)\\
   	\dfrac{1}{f(t)} & 0\\
\end{matrix}  & 0  \\
\hline
0 & 0 & I_{n-i-1}
\end{array} \right) \text{ for } 1\leq i \leq n-1,$$
where $f(t)$ is invertible in $ \mathbb{Z}[t^{\pm 1}]$.
\end{definition}

\begin{theorem} \cite{Nasser2025}
The representation $\eta_2: T_n \longrightarrow \mathrm{GL}_n(\mathbb{Z}[t^{\pm 1}])$ is reducible. 
\end{theorem}

\begin{theorem} \cite{Nasser2025}
The representation $\eta_2: T_n \longrightarrow \mathrm{GL}_n(\mathbb{Z}[t^{\pm 1}])$ is faithful for $n=2$ and unfaithful for all $n\geq 3$.
\end{theorem}

\section{On the Composition of the $\eta_1$ Representation and the Cohen Map}
 
In this section, we construct a high-dimensional linear representation of $T_n$ induced from the representation $\eta_1$ given in Definition \ref{defeta1}. The construction relies on the application of the Cohen map to the twin group $T_n$.

\begin{definition} \cite{Cohen1996} \label{defcoh}
The Cohen map is a homomorphism $\mathcal{C}: T_n \longrightarrow T_{nk}$ defined by
\vskip -12pt
\[
s_i \longmapsto r_i
=
(s_{ki}s_{ki+1}\cdots s_{ki+k-1})
(s_{ki-1}s_{ki}\cdots s_{ki+k-2})\cdots
(s_{ki-k+1}s_{ki-k+2}\cdots s_{ki}),
\]
where $k$ is a positive integer. The twin $r_i$ can be interpreted geometrically as the twin obtained by replacing each strand of the geometric twin $s_i$ with $k$ parallel strands. 
\end{definition}

In this paper, we focus on the case $k=2$, aiming to construct a linear representation of $T_n$ of degree $2n$. That is, we take $$r_i=s_{2i}s_{2i+1}s_{2i-1}s_{2i}.$$ Let $x_1, x_2, \dots, x_{2n}$ be a fixed generating set of the free group $\mathbb{F}_{2n}$. We consider the composition
\[
T_n \xrightarrow{\ \mathcal{C}\ } T_{2n}
\xrightarrow{\ \eta_1 \ } \mathrm{Aut}(\mathbb{F}_{2n})
\xrightarrow{\ \mathcal{M}\ } \mathrm{GL}_{2n}(\mathbb{Z}[t^{\pm1}]).
\]
where
\begin{enumerate}
    \item the first map is the Cohen map $\mathcal{C}$ of $T_n$ (Definition \ref{defcoh}),
    \item the second map is the representation $\eta_1$ of $T_{2n}$ (Definition \ref{defeta1}),
\item The third map is the Magnus representation, denoted by $\mathcal{M}$, which assigns to each automorphism in $\mathrm{Aut}(\mathbb{F}_{2n})$ a $2n \times 2n$ invertible matrix over the Laurent polynomial ring $\mathbb{Z}[t^{\pm1}]$, using Fox derivatives. Here, $t$ is an indeterminate corresponding to the image of the generators $y_1, y_2, \ldots, y_{2n}$ of $\mathbb{F}_{2n}$ under the map $\Psi:\mathbb{F}_{2n} \longrightarrow \mathbb{Z}[t^{\pm1}]$, that is, $\Psi(y_i) = t$ for all $1 \le i \le 2n$. Hence, each generator $y_i$ is replaced by $t$ in the matrix entries.
\end{enumerate}

\noindent Our procedure can be described in the following three steps. We consider each step and its detailed calculations and proofs as follows. \vspace*{0.1cm}

\noindent \underline{Step 1.} We choose a new generating set for the free group $\mathbb{F}_{2n}$. Consider the generating set $y_1,y_2,\ldots, y_{2n}$, where, for $1\leq i \leq n$,
\[
y_i = x_{2i-1} \  \text{ and } \  y_{n+i} = x_{2i}.
\]

\noindent \underline{Step 2.} We compute the action of each automorphism associated with $r_i=\mathcal{C}(s_i)$ on the chosen generating set $y_1, y_2, \ldots, y_{2n}$. The result is in the following proposition.

\begin{proposition} \label{propppp}
For $1\leq i \leq n-1$ and $1\leq j\leq 2n$, the action of $r_i$ on the generators $y_j$ of $\mathbb{F}_{2n}$ described in Step 1, under the representation $\eta_1$ introduced in Definition \ref{defeta1}, is given as follows.
\begin{itemize}
\item[(1)] $y_i \longmapsto y_iy_{n+i}y_{i+1}y_{n+i}^{-1}y_i^{-1}.$
\item[(2)] $y_{i+1} \longmapsto y_{i}y_{n+i}y_{i+1}^{-1}y_{n+i}^{-1}y_{i}y_{n+i}y_{i+1}y_{n+i}^{-1}y_{i}^{-1}y_{n+i}y_{i+1}^{-1}y_{n+i}y_{i+1}y_{n+i}^{-1}y_{n+i+1}^{-1}y_{n+i}\\ \hspace*{1.4cm} y_{i+1}^{-1}y_{n+i}^{-1}y_{i+1}y_{n+i}^{-1}y_{i}y_{n+i}y_{i+1}^{-1} y_{n+i}^{-1}y_{i}y_{n+i}y_{i+1}y_{n+i}^{-1}y_{i}^{-1}y_{n+i}y_{i+1}^{-1}y_{n+i}y_{i+1}\\ \hspace*{1.4cm} y_{n+i}^{-1}y_{n+i+1}y_{n+i}y_{i+1}^{-1}y_{n+i}^{-1}y_{i+1}y_{n+i}^{-1}y_{i}y_{n+i}y_{i+1}^{-1}y_{n+i}^{-1}y_{i}^{-1}y_{n+i}y_{i+1}y_{n+i}^{-1}y_{i}^{-1}.$
\item[(3)] $y_{n+i} \longmapsto y_{i}y_{n+i}y_{i+1}^{-1}y_{n+i}^{-1}y_{i}y_{n+i}y_{i+1}y_{n+i}^{-1}y_{i}^{-1}y_{n+i}y_{i+1}^{-1}y_{n+i}y_{i+1}y_{n+i}^{-1}y_{n+i+1}y_{n+i}\\ \hspace*{1.4cm} y_{i+1}^{-1}y_{n+i}^{-1}y_{i+1}y_{n+i}^{-1}y_{i}y_{n+i}y_{i+1}^{-1}y_{n+i}^{-1}y_{i}^{-1}y_{n+i}y_{i+1}y_{n+i}^{-1}y_{i}^{-1}.  $
\item[(4)] $y_{n+i+1} \longmapsto y_{n+i}y_{i+1}^{-1}y_{n+i}y_{i+1}y_{n+i}^{-1}y_{n+i+1}^{-1}y_{n+i}y_{i+1}^{-1}y_{n+i}y_{i+1}y_{n+i}^{-1}y_{n+i+1}y_{n+i}y_{i+1}^{-1}\\ \hspace*{1.775cm} y_{n+i}^{-1}y_{i+1}y_{n+i}^{-1}.$
\item[(5)] $y_{\ell} \longmapsto y_\ell$ for all $\ell \notin \{i,i+1,n+i,n+i+1\}$. 
\end{itemize}

\begin{proof}
Fix $1\leq i \leq n-1$. We consider each case separately in the following.
\begin{itemize}
\item[(1)] $r_i(y_i)=s_{2i}s_{2i+1}s_{2i-1}s_{2i}(y_i)\\
\hspace*{0.83cm} =s_{2i}s_{2i+1}s_{2i-1}s_{2i}(x_{2i-1})\\
\hspace*{0.83cm} =s_{2i}s_{2i+1}s_{2i-1}(x_{2i-1})\\
\hspace*{0.83cm} =s_{2i}s_{2i+1}(x_{2i-1}x_{2i}x_{2i-1}^{-1})\\
\hspace*{0.83cm} =s_{2i}(x_{2i-1}x_{2i}x_{2i-1}^{-1})\\
\hspace*{0.83cm} =x_{2i-1}x_{2i}x_{2i+1}x_{2i}^{-1}x_{2i-1}^{-1}\\
\hspace*{0.83cm} =y_iy_{n+i}y_{i+1}y_{n+i}^{-1}y_i^{-1}.\vspace*{0.1cm}
$
\item[(2)] $r_i(y_{i+1})=s_{2i}s_{2i+1}s_{2i-1}s_{2i}(y_{i+1})\\
\hspace*{1.2cm}=s_{2i}s_{2i+1}s_{2i-1}s_{2i}(x_{2i+1})\\
\hspace*{1.2cm}=s_{2i}s_{2i+1}s_{2i-1}(x_{2i}x_{2i+1}^{-1}x_{2i}x_{2i+1}x_{2i}^{-1})\\
\hspace*{1.2cm}=s_{2i}s_{2i+1}(x_{2i-1}x_{2i}^{-1}x_{2i-1}x_{2i}x_{2i-1}^{-1}x_{2i+1}^{-1}x_{2i-1}x_{2i}^{-1}x_{2i-1}x_{2i}x_{2i-1}^{-1}x_{2i+1}\\
\hspace*{3cm} x_{2i-1}x_{2i}^{-1}x_{2i-1}^{-1}x_{2i}x_{2i-1}^{-1})\\
\hspace*{1.2cm}=s_{2i}(x_{2i-1}x_{2i}^{-1}x_{2i-1}x_{2i}x_{2i-1}^{-1}x_{2i+1}x_{2i+2}^{-1}x_{2i+1}^{-1}x_{2i-1}x_{2i}^{-1}x_{2i-1}x_{2i}x_{2i-1}^{-1}\\
\hspace*{2.22cm} x_{2i+1}x_{2i+2}x_{2i+1}^{-1}x_{2i-1}x_{2i}^{-1}x_{2i-1}^{-1}x_{2i}x_{2i-1}^{-1})\\
\hspace*{1.2cm}=x_{2i-1}x_{2i}x_{2i+1}^{-1}x_{2i}^{-1}x_{2i-1}x_{2i}x_{2i+1}x_{2i}^{-1}x_{2i-1}^{-1}x_{2i}x_{2i+1}^{-1}x_{2i}x_{2i+1}x_{2i}^{-1}x_{2i+2}^{-1}\\
\hspace*{1.65cm} x_{2i}x_{2i+1}^{-1}x_{2i}^{-1}x_{2i+1}x_{2i}^{-1}x_{2i-1}x_{2i}x_{2i+1}^{-1}x_{2i}^{-1}x_{2i-1}x_{2i}x_{2i+1}x_{2i}^{-1}x_{2i-1}^{-1}x_{2i}\\
\hspace*{1.65cm} x_{2i+1}^{-1}x_{2i}x_{2i+1}x_{2i}^{-1}x_{2i+2}x_{2i}x_{2i+1}^{-1}x_{2i}^{-1}x_{2i+1}x_{2i}^{-1}x_{2i-1}x_{2i}x_{2i+1}^{-1}x_{2i}^{-1}x_{2i-1}^{-1}\\
\hspace*{1.65cm} x_{2i}x_{2i+1}x_{2i}^{-1}x_{2i-1}^{-1}\\
\hspace*{1.2cm}=y_{i}y_{n+i}y_{i+1}^{-1}y_{n+i}^{-1}y_{i}y_{n+i}y_{i+1}y_{n+i}^{-1}y_{i}^{-1}y_{n+i}y_{i+1}^{-1}y_{n+i}y_{i+1}y_{n+i}^{-1}y_{n+i+1}^{-1}y_{n+i}\\ \hspace*{1.67cm} y_{i+1}^{-1}y_{n+i}^{-1}y_{i+1}y_{n+i}^{-1}y_{i}y_{n+i}y_{i+1}^{-1} y_{n+i}^{-1}y_{i}y_{n+i}y_{i+1}y_{n+i}^{-1}y_{i}^{-1}y_{n+i}y_{i+1}^{-1}y_{n+i}y_{i+1}\\ \hspace*{1.67cm} y_{n+i}^{-1}y_{n+i+1}y_{n+i}y_{i+1}^{-1}y_{n+i}^{-1}y_{i+1}y_{n+i}^{-1}y_{i}y_{n+i}y_{i+1}^{-1}y_{n+i}^{-1}y_{i}^{-1}y_{n+i}y_{i+1}y_{n+i}^{-1}y_{i}^{-1}.\vspace*{0.1cm}
$
\item[(3)] $r_i(y_{n+i})=s_{2i}s_{2i+1}s_{2i-1}s_{2i}(y_{n+i})\\
\hspace*{1.225cm}=s_{2i}s_{2i+1}s_{2i-1}s_{2i}(x_{2i})\\
\hspace*{1.225cm}=s_{2i}s_{2i+1}s_{2i-1}(x_{2i}x_{2i+1}x_{2i}^{-1})\\
\hspace*{1.225cm}=s_{2i}s_{2i+1}(x_{2i-1}x_{2i}^{-1}x_{2i-1}x_{2i}x_{2i-1}^{-1}x_{2i+1}x_{2i-1}x_{2i}^{-1}x_{2i-1}^{-1}x_{2i}x_{2i-1}^{-1})\\
\hspace*{1.225cm}=s_{2i}(x_{2i-1}x_{2i}^{-1}x_{2i-1}x_{2i}x_{2i-1}^{-1}x_{2i+1}x_{2i+2}x_{2i+1}^{-1}x_{2i-1}x_{2i}^{-1}x_{2i-1}^{-1}x_{2i}x_{2i-1}^{-1})\\
\hspace*{1.22cm}=x_{2i-1}x_{2i}x_{2i+1}^{-1}x_{2i}^{-1}x_{2i-1}x_{2i}x_{2i+1}x_{2i}^{-1}x_{2i-1}^{-1}x_{2i}x_{2i+1}^{-1}x_{2i}x_{2i+1}x_{2i}^{-1}x_{2i+2}\\
\hspace*{1.68cm}x_{2i}x_{2i+1}^{-1}x_{2i}^{-1}x_{2i+1}x_{2i}^{-1}x_{2i-1}x_{2i}x_{2i+1}^{-1}x_{2i}^{-1}x_{2i-1}^{-1}x_{2i}x_{2i+1}x_{2i}^{-1}x_{2i-1}^{-1}\\
\hspace*{1.2cm} =y_{i}y_{n+i}y_{i+1}^{-1}y_{n+i}^{-1}y_{i}y_{n+i}y_{i+1}y_{n+i}^{-1}y_{i}^{-1}y_{n+i}y_{i+1}^{-1}y_{n+i}y_{i+1}y_{n+i}^{-1}y_{n+i+1}y_{n+i}\\
\hspace*{1.66cm} y_{i+1}^{-1}y_{n+i}^{-1}y_{i+1}y_{n+i}^{-1}y_{i}y_{n+i}y_{i+1}^{-1}y_{n+i}^{-1}y_{i}^{-1}y_{n+i}y_{i+1}y_{n+i}^{-1}y_{i}^{-1}.\vspace*{0.1cm}
$
\item[(4)] $r_i(y_{n+i+1})=s_{2i}s_{2i+1}s_{2i-1}s_{2i}(y_{n+i+1})\\
\hspace*{1.575cm}=s_{2i}s_{2i+1}s_{2i-1}s_{2i}(x_{2i+2})\\
\hspace*{1.575cm}=s_{2i}s_{2i+1}s_{2i-1}(x_{2i+2})\\
\hspace*{1.575cm}=s_{2i}s_{2i+1}(x_{2i+2})\\
\hspace*{1.575cm}=s_{2i}(x_{2i+1}x_{2i+2}^{-1}x_{2i+1}x_{2i+2}x_{2i+1}^{-1})\\
\hspace*{1.575cm}=x_{2i}x_{2i+1}^{-1}x_{2i}x_{2i+1}x_{2i}^{-1}x_{2i+2}^{-1}x_{2i}x_{2i+1}^{-1}x_{2i}x_{2i+1}x_{2i}^{-1}x_{2i+2}x_{2i}x_{2i+1}^{-1}x_{2i}^{-1}\\
\hspace*{2.05cm} x_{2i+1}x_{2i}^{-1}\\
\hspace*{1.575cm}= y_{n+i}y_{i+1}^{-1}y_{n+i}y_{i+1}y_{n+i}^{-1}y_{n+i+1}^{-1}y_{n+i}y_{i+1}^{-1}y_{n+i}y_{i+1}y_{n+i}^{-1}y_{n+i+1}y_{n+i}y_{i+1}^{-1}\\ \hspace*{2.05cm} y_{n+i}^{-1}y_{i+1}y_{n+i}^{-1}.$\vspace*{0.1cm}
\item[(5)] $r_i(y_{\ell})=y_{\ell}$ for every $\ell \notin \{i,i+1,n+i,n+i+1\}$ follows directly, since $r_i$  only acts on $y_i$, $y_{i+1}$, $y_{n+i}$, and $y_{n+i+1}$. This is because the generator $s_i$ acts only on $x_i$ and $x_{i+1}$, leaving all other generators unchanged.
\end{itemize}
\end{proof}
\end{proposition}

\noindent \underline{Step 3.} Apply Fox derivatives to the resulting automorphisms in $\mathrm{Aut}(\mathbb{F}_{2n})$ obtained in Proposition \ref{propppp}, which produces matrices over $\mathbb{Z}[t^{\pm1}]$. We start by the definition of Fox derivatives.

\begin{definition}
Let $\mathbb{F}_n$ be the free group with $n$ generators 
$x_1, x_2, \ldots, x_n$. For $j=1,2,\ldots,n$, the Fox derivatives on the group ring $\mathbb{Z}[\mathbb{F}_n]$ are defined as follows.
\begin{itemize}
\item[(1)] $\displaystyle \frac{\partial}{\partial x_j}x_i=\delta_{ij}$,
\item[(2)] $\displaystyle \frac{\partial }{\partial x_j}x_i^{-1}=-\delta_{ij}x_i^{-1}$,
\item[(3)] $\displaystyle \frac{\partial}{\partial x_j}(uv)=\frac{\partial }{\partial x_j}u\,\epsilon(v)+u\frac{\partial }{\partial x_j}v$ for $u,v\in \mathbb{Z}[\mathbb{F}_n]$.
\end{itemize}
Here, $\delta_{ij}$ denotes the Kronecker symbol and $\epsilon(v)=1$ whenever $v\in \mathbb{F}_n$.
\end{definition}

Now, we set $D_\ell=\Psi \displaystyle \frac{\partial}{\partial y_\ell}$. Recall that $\Psi$ is defined by $\Psi(y_i)=t$ for all $1\leq i \leq 2n$. For each $r_i$, $1\leq i \leq n-1$, we compute its corresponding Jacobian matrix as follows:
\[
J(r_i)=
\begin{pmatrix}
D_1\big(r_i(y_1)\big) & \cdots & D_{2n}\big(r_i(y_1)\big) \\
\vdots & & \vdots \\
D_1\big(r_i(y_{2n})\big) & \cdots & D_{2n}\big(r_i(y_{2n})\big)
\end{pmatrix}.
\]
We now prove our main theorem in this section.

\begin{theorem} \label{Th1st}
The linear representation of $T_n$ obtained by composing the Cohen map given in Definition \ref{defcoh} with the $\eta_1$ representation given in Definition \ref{defeta1} is the representation 
$$\widetilde{\eta}_1:T_n \longrightarrow \mathrm{GL}_{2n}(\mathbb{Z}[t^{\pm1}])$$ that acts on the generators $s_i, 1\leq i \leq n-1$, as follows.
$$\widetilde{\eta}_1(s_i)=$$
{\small
\[
\left(
\begin{array}{c|c|c}
I_{i-1} & 0 & 0 \\
\hline
0 &
\begin{matrix}
1-t & t^2 & 0 & \cdots & 0 & t(1-t) & 0\\
(t-2)^2 & -2t^3+6t^2-5t+1 & 0 & \cdots & 0 & 2t^3-8t^2+10t-4 & t(t-1)\\
0 & 0 & 1 & \cdots & 0 & 0 & 0\\
\vdots & \vdots & & \ddots & & \vdots & \vdots\\
0 & 0 & 0 & \cdots & 1 & 0 & 0\\
(t-1)(t-2) & -2t(t-1)^2 & 0 & \cdots & 0 & 2t^3-6t^2+5t-1 & t^2\\
0 & (1-t)(t-2) & 0 & \cdots & 0 & (t-2)^2 & t-1
\end{matrix}
& 0 \\
\hline
0 & 0 & I_{n-i-1}
\end{array}
\right).
\]
}
\end{theorem} 

\begin{proof}
Fix $1\le i \le n-1$. First, we compute 
\[
\frac{\partial }{\partial y_\ell}r_i(y_j)
\]
for all $1\le j,\ell \le 2n$. Next, we obtain $D_{\ell}( r_i(y_j))$ for all $1\le j,\ell \le 2n$ by replacing each $y_j$ by the indeterminate $t$. For every $1\leq j \leq 2n$, we have to consider the following five cases:
$$\ell=i,\ \ell=i+1,\ \ell=n+i,\ \ell=n+i+1, \ \ell \notin \{i, i+1, n+i, n+i+1\}.$$
Since these computations are tough, we illustrate the computation only for the case $j=i$; that is, we compute the $i^{\text{th}}$ row of $\widetilde{\eta}_1(s_i)$, while the other rows can be obtained similarly. Note that a straightforward verification shows that this matrix satisfies both families of relations that define $T_n$, confirming that it is a well-constructed representation. So, we now fix $j=i$ and we consider the following five cases:
\begin{itemize}
\item[(1)] $\displaystyle \frac{\partial }{\partial y_i}r_i(y_i)=\frac{\partial }{\partial y_i}(y_iy_{n+i}y_{i+1}y_{n+i}^{-1}y_i^{-1})\\
\hspace*{1.41cm} =\frac{\partial }{\partial y_i}(y_i)+y_i\frac{\partial }{\partial y_i}(y_{n+i}y_{i+1}y_{n+i}^{-1}y_i^{-1})\\
\hspace*{1.41cm} =1+y_iy_{n+i}y_{i+1}y_{n+i}^{-1}\frac{\partial }{\partial y_i}(y_i^{-1})\\
\hspace*{1.41cm} =1-y_iy_{n+i}y_{i+1}y_{n+i}^{-1}y_i^{-1}$.

\noindent Now, replacing each $y_{m}$ by $t$, we get that $D_{i}( r_i(y_i))=1-t$. Thus, the $i\times i$ entry in the matrix $\widetilde{\eta}_1(s_i)$ is $1-t$, as required.
\vspace*{0.1cm}
\item[(2)] $ \displaystyle \frac{\partial }{\partial y_{i+1}}r_i(y_{i})=\frac{\partial }{\partial y_{i+1}}(y_iy_{n+i}y_{i+1}y_{n+i}^{-1}y_i^{-1})\\
\hspace*{1.77cm} =y_iy_{n+i}\frac{\partial }{\partial y_{i+1}}(y_{i+1}y_{n+i}^{-1}y_i^{-1})\\
\hspace*{1.77cm} =y_iy_{n+i}\left(\frac{\partial }{\partial y_{i+1}}y_{i+1} + y_{i+1}\frac{\partial }{\partial y_{i+1}}(y_{n+i}^{-1}y_i^{-1})\right)\\
\hspace*{1.77cm} =y_iy_{n+i}(1+0)\\
\hspace*{1.77cm}=y_iy_{n+i}.$

\noindent Now, replacing each $y_{m}$ by $t$, we get that $D_{i+1}( r_i(y_i))=t^2$. Thus, the $i\times (i+1)$ entry in the matrix $\widetilde{\eta}_1(s_i)$ is $t^2$, as required.
\vspace*{0.1cm}
\item[(3)] $\displaystyle \frac{\partial }{\partial y_{n+i}}r_i(y_i)=\frac{\partial }{\partial y_{n+i}}(y_iy_{n+i}y_{i+1}y_{n+i}^{-1}y_i^{-1})\\
\hspace*{1.77cm}=y_i\frac{\partial }{\partial y_{n+i}}(y_{n+i}y_{i+1}y_{n+i}^{-1}y_i^{-1})\\
\hspace*{1.77cm}=y_i\left( \frac{\partial }{\partial y_{n+i}}y_{n+i}+y_{n+i}\frac{\partial }{\partial y_{n+i}}(y_{i+1}y_{n+i}^{-1}y_i^{-1})\right)\\
\hspace*{1.77cm}=y_i\left( 1+y_{n+i}y_{i+1}\frac{\partial }{\partial y_{n+i}}(y_{n+i}^{-1}y_i^{-1})\right)\\
\hspace*{1.77cm}=y_i \left( 1+y_{n+i}y_{i+1}\left( \frac{\partial }{\partial y_{n+i}}y_{n+i}^{-1}+y_{n+i}^{-1}\frac{\partial }{\partial y_{n+i}}y_i^{-1}\right) \right)\\
\hspace*{1.77cm}=y_i\left( 1+y_{n+i}y_{i+1}(-y_{n+i}^{-1}+0)\right)\\
\hspace*{1.77cm}=y_i( 1-y_{n+i}y_{i+1}y_{n+i}^{-1}).\\
$
Now, replacing each $y_{m}$ by $t$, we get that $D_{n+i}( r_i(y_i))=t(1-t)$. Thus, the $i\times (n+i)$ entry in the matrix $\widetilde{\eta}_1(s_i)$ is $t(1-t)$, as required.
\vspace*{0.1cm}
\item[(4)] $\displaystyle  \frac{\partial }{\partial y_{n+i+1}}r_i(y_i)=\frac{\partial }{\partial y_{n+i+1}}(y_iy_{n+i}y_{i+1}y_{n+i}^{-1}y_i^{-1})=0.$\\
Thus, we get that $D_{n+i+1}( r_i(y_i))=0$, and so the $i\times (n+i+1)$ entry in the matrix $\widetilde{\eta}_1(s_i)$ is $0$, as required.
\item[(5)] $\displaystyle  \frac{\partial }{\partial y_{\ell}}r_i(y_i)=\frac{\partial}{\partial y_{\ell}}(y_iy_{n+i}y_{i+1}y_{n+i}^{-1}y_i^{-1})=0.$ for all $\ell \notin \{i, i+1, n+i, n+i+1\}$.\\
Thus, we get that $D_{\ell}( r_i(y_i))=0$, and so the $i\times \ell$ entry in the matrix $\widetilde{\eta}_1(s_i)$ is $0$ for all $\ell \notin \{i, i+1, n+i, n+i+1\}$, as required.
\end{itemize}
A similar computation for all $2n$ rows can be carried out in the same manner, and thus the proof is completed.
\end{proof}

\section{On the Tensoring of the $\eta_2$ Representation}

In this section, we construct a high-dimensional representation of $T_n$ based on the representation $\eta_2$ introduced in Definition \ref{defeta2}, using the tensor product of representations. Here, we focus on the cases $n=2$ and $n=3$, although the same approach can be applied to larger values of $n$.

\begin{definition}
Let $\eta: G\to \mathrm{GL}(V)$ and $\eta': G\to \mathrm{GL}(W)$ be two representations of a group $G$. The tensor product representation
\[
\eta\otimes\eta' : G\to \mathrm{GL}(V\otimes W)
\]
is defined by
\[
(\eta\otimes\eta')(g)(v\otimes w)=\eta(g)v\otimes\eta'(g)w.
\]
In matrix form, the tensor product representation is obtained by taking the Kronecker product of the matrices corresponding to $\eta(g)$ and $\eta'(g)$.
\end{definition}

Remark that if $A=(a_{ij})$ is a matrix of size $m\times n$ and $B$ is a matrix of size $p\times q$, the Kronecker product $A\otimes B$ is the block matrix obtained by multiplying each entry $a_{ij}$ of $A$ by the matrix $B$. More precisely,
\[
A\otimes B=
\begin{pmatrix}
a_{11}B & a_{12}B & \cdots & a_{1n}B\\
a_{21}B & a_{22}B & \cdots & a_{2n}B\\
\vdots & \vdots & \ddots & \vdots\\
a_{m1}B & a_{m2}B & \cdots & a_{mn}B
\end{pmatrix}.
\]

We now compute the tensor product of the representation $\eta_2$ with itself for $n=2$.

\begin{proposition} \label{firstprop}
Set $n=2$ and consider the representation $\eta_2:T_2 \longrightarrow \mathrm{GL}_2(\mathbb{Z}[t^{\pm 1}])$ given in Definition \ref{defeta2}. The tensor product of $\eta_2$ with itself is the representation $$\hat{\eta}_2:T_2 \longrightarrow \mathrm{GL}_4(\mathbb{Z}[t^{\pm 1}])$$ given by
$$
\hat{\eta}_2(s_1)=
\left( \begin{matrix}
0 & 0 & 0 & f^2(t)\\
0 & 0 & 1 & 0\\
0 & 1 & 0 & 0\\
\dfrac{1}{f^2(t)} & 0 & 0 & 0\\
\end{matrix}\right),$$
where $f(t)$ is invertible in $ \mathbb{Z}[t^{\pm 1}]$.
\end{proposition}

\begin{proof}
From Definition \ref{defeta2}, we have
$$
\eta_2(s_1)=
\left( \begin{matrix}
0 & f(t)\\
\dfrac{1}{f(t)} & 0
\end{matrix}\right) ,
$$
where $f(t)$ is invertible in $ \mathbb{Z}[t^{\pm 1}]$. So, we obtain that
$$
\hat{\eta}_2(s_1)=
\left( \begin{matrix}
0\ \eta_2(s_1) & f(t)\ \eta_2(s_1)\\
\dfrac{1}{f(t)}\ \eta_2(s_1) & 0\ \eta_2(s_1)
\end{matrix}\right)=\left( \begin{matrix}
0\ \left( \begin{matrix}
0 & f(t)\\
\dfrac{1}{f(t)} & 0
\end{matrix}\right) & f(t)\ \left( \begin{matrix}
0 & f(t)\\
\dfrac{1}{f(t)} & 0
\end{matrix}\right)\\
\dfrac{1}{f(t)}\ \left( \begin{matrix}
0 & f(t)\\
\dfrac{1}{f(t)} & 0
\end{matrix}\right) & 0\ \left( \begin{matrix}
0 & f(t)\\
\dfrac{1}{f(t)} & 0
\end{matrix}\right)
\end{matrix}\right).
$$
Thus, we get the desired matrix
$$
\hat{\eta}_2(s_1)=
\left( \begin{matrix}
0 & 0 & 0 & f^2(t)\\
0 & 0 & 1 & 0\\
0 & 1 & 0 & 0\\
\dfrac{1}{f^2(t)} & 0 & 0 & 0\\
\end{matrix}\right).
$$
\end{proof}

\begin{proposition}
The representation $\hat{\eta}_2:T_2 \longrightarrow \mathrm{GL}_4(\mathbb{Z}[t^{\pm 1}])$ is faithful and reducible.
\end{proposition}

\begin{proof}
It is clear that the representation is faithful since $T_2=\langle s_1 \rangle$ and $\hat{\eta}_2(s_1)\neq I_4$, the $4\times 4$ identity matrix. Moreover, the two subspaces $U_1=\langle e_1, e_4 \rangle$ and $U_2=\langle e_2, e_3 \rangle$ of $\mathbb{C}^4$ are invariant under the action of $\hat{\eta}_2(s_1)$, where $e_i$ denote the standard basis vectors of $\mathbb{C}^4$. Therefore, the representation $\hat{\eta}_2$ is reducible.
\end{proof}

We now compute the tensor product of the representation $\eta_2$ with itself for $n=3$.

\begin{proposition} \label{secondprop}
Set $n=3$ and consider the representation $\eta_2:T_3 \longrightarrow \mathrm{GL}_3(\mathbb{Z}[t^{\pm 1}])$ given in Definition \ref{defeta2}. The tensor product of $\eta_2$ to itself is the representation $$\bar{\eta}_2:T_3 \longrightarrow \mathrm{GL}_9(\mathbb{Z}[t^{\pm 1}])$$ such that
$$
\bar{\eta}_2(s_1)=
\left( \begin{matrix}
0 & 0 & 0 & 0 & f^2(t) & 0 & 0 & 0 & 0\\
0 & 0 & 0 & 1 & 0 & 0 & 0 & 0 & 0\\
0 & 0 & 0 & 0 & 0 & f(t) & 0 & 0 & 0\\
0 & 1 & 0 & 0 & 0 & 0 & 0 & 0 & 0\\
\dfrac{1}{f^2(t)} & 0 & 0 & 0 & 0 & 0 & 0 & 0 & 0\\
0 & 0 & \dfrac{1}{f(t)} & 0 & 0 & 0 & 0 & 0 & 0\\
0 & 0 & 0 & 0 & 0 & 0 & 0 & f(t) & 0\\
0 & 0 & 0 & 0 & 0 & 0 & \dfrac{1}{f(t)} & 0 & 0\\
0 & 0 & 0 & 0 & 0 & 0 & 0 & 0 & 1\\
\end{matrix}\right)
$$
and 
$$
\bar{\eta}_2(s_2)=
\left( \begin{matrix}
1 & 0 & 0 & 0 & 0 & 0 & 0 & 0 & 0\\
0 & 0 & f(t) & 0 & 0 & 0 & 0 & 0 & 0\\
0 & \dfrac{1}{f(t)} & 0 & 0 & 0 & 0 & 0 & 0 & 0\\
0 & 0 & 0 & 0 & 0 & 0 & f(t) & 0 & 0\\
0 & 0 & 0 & 0 & 0 & 0 & 0 & 0 & f^2(t)\\
0 & 0 & 0 & 0 & 0 & 0 & 0 & 1 & 0\\
0 & 0 & 0 & \dfrac{1}{f(t)} & 0 & 0 & 0 & 0 & 0\\
0 & 0 & 0 & 0 & 0 & 1 & 0 & 0 & 0\\
0 & 0 & 0 & 0 & \dfrac{1}{f^2(t)} & 0 & 0 & 0 & 0\\
\end{matrix}\right),
$$
where $f(t)$ is invertible in $ \mathbb{Z}[t^{\pm 1}]$.
\end{proposition}

\begin{proof}
The argument is analogous to that in Proposition \ref{firstprop} and is thus omitted.
\end{proof}

\begin{proposition}
The representation $\bar{\eta}_2:T_3 \longrightarrow \mathrm{GL}_9(\mathbb{Z}[t^{\pm 1}])$ is unfaithful and reducible.
\end{proposition}

\begin{proof}
The representation $\bar{\eta}_2$ is unfaithful because $\bar{\eta}_2(s_1s_2s_1s_2s_1s_2) = I_9$, while $s_1s_2s_1s_2s_1s_2$ is a nontrivial element in $T_3$, since the relation $s_i s_{i+1} s_i = s_{i+1} s_i s_{i+1}$ does not hold in $T_n$. Furthermore, the subspace $U = \langle e_1, e_5, e_9 \rangle$ of $\mathbb{C}^9$ is invariant under $\bar{\eta}_2(s_1)$ and $\bar{\eta}_2(s_2)$, where $e_i$ denote the standard basis vectors of $\mathbb{C}^9$. Therefore, the representation $\bar{\eta}_2$ is reducible.
\end{proof}

\section{Future Work}

This section presents directions for research inspired by the techniques introduced in this paper.
\begin{itemize}
    \item[(F1)] \underline{ The Cohen map:} Explore larger values of $k$ in the Cohen map to construct higher-dimensional representations in a similar way of constructing $\widetilde{\eta}_1$.\vspace{0.1cm}
    \item[(F2)] \underline{Detailed study of $\widetilde{\eta}_1$:} Explore structural properties of the new representation $\widetilde{\eta}_1$ of $T_n$, including its faithfulness and reducibility.\vspace{0.1cm}
    \item[(F3)] \underline{The tensor product of $\eta_2$ for larger \(n\):} Explore the tensor product of the representation $\eta_2$ of $T_n$ with itself for $n>3$ and study the faithfulness and reducibility.
    \vspace{0.1cm}
    \item[(F4)] \underline{Application to other groups:} Utilize these methods to construct representations for algebraic structures beyond $T_n$, for instance the triplet group $L_n$, and investigate their properties. Representations of $L_n$ can be found in sources such as \cite{PKumar2024,Nasser2026}.
\end{itemize}

\end{document}